\documentclass[12pt]{amsart}
\usepackage{amscd}
\usepackage[latin1]{inputenc}
\usepackage{amsmath}
\usepackage{amssymb}
\usepackage{amsthm}
\usepackage{latexsym}
\usepackage{pstcol,pst-plot,pst-3d}
\usepackage{mathptmx}
\usepackage{multicol}
\usepackage{hyperref}

\psset{unit=0.7cm,linewidth=0.8pt,arrowsize=2.5pt 4}

\newpsstyle{fatline}{linewidth=1.5pt}
\newpsstyle{fyp}{fillstyle=solid,fillcolor=verylight}
\definecolor{verylight}{gray}{0.97}
\definecolor{light}{gray}{0.9}
\definecolor{medium}{gray}{0.85}

\def\NZQ{\Bbb}               
\def\NN{{\NZQ N}}

\def\PP{{\NZQ P}}
\def\frk{\frak}               

\def\mm{{\frk m}}

\def\Phi{{\frk n}}
\def\Phi{{\frk N}}
\def\opn#1#2{\def#1{\operatorname{#2}}} 
\opn\chara{char} \opn\length{\ell} \opn\pd{pd} \opn\rk{rk}
\opn\projdim{proj\,dim} \opn\injdim{inj\,dim} \opn\rank{rank}
\opn\depth{depth} \opn\grade{grade} \opn\height{height}
\opn\embdim{emb\,dim} \opn\codim{codim} \opn\dim{dim}

\opn\Tr{Tr} \opn\bigrank{big\,rank}
\opn\superheight{superheight}\opn\lcm{lcm}
\opn\trdeg{tr\,deg}%
\opn\reg{reg} \opn\lreg{lreg} \opn\ini{in}
\opn\div{div} \opn\Div{Div} \opn\cl{cl} \opn\Cl{Cl}
\opn\Spec{Spec} \opn\Supp{Supp} \opn\supp{supp} \opn\Sing{Sing}
\opn\Ass{Ass}  \opn\Min{Min}
\opn\Ann{Ann} \opn\Rad{Rad} \opn\Soc{Soc}
\opn\Ker{Ker} \opn\Coker{Coker} \opn\Am{Am} \opn\Hom{Hom}
\opn\Tor{Tor} \opn\Ext{Ext} \opn\End{End} \opn\Aut{Aut}
\opn\id{id}  \opn\e{e}

\opn\nat{nat} \opn\deg{deg} \opn\adeg{adeg} \opn\Hilb{Hilb}
\opn\pff{pf}
\opn\Pf{Pf} \opn\GL{GL} \opn\SL{SL} \opn\mod{mod} \opn\ord{ord}
\opn\aff{aff} \opn\con{conv} \opn\relint{relint} \opn\st{st}
\opn\lk{lk} \opn\cn{cn} \opn\core{core} \opn\vol{vol}
\opn\link{link}  \opn\infpt{infpt} \opn\mult{mult}
\opn\gr{gr} \opn\Std{Std}

\def\pot#1#2{#1[\kern-0.28ex[#2]\kern-0.28ex]}

\opn\dirlim{\underrightarrow{\lim}}
\opn\inivlim{\underleftarrow{\lim}}
\let\union=\cup
\let\sect=\cap
\let\dirsum=\oplus

\let\iso=\cong
\let\Union=\bigcup
\let\Sect=\bigcap
\let\Dirsum=\bigoplus

\let\to=\rightarrow

\def\Implies{\ifmmode\Longrightarrow \else
     \unskip${}\Longrightarrow{}$\ignorespaces\fi}
\def\implies{\ifmmode\Rightarrow \else
     \unskip${}\Rightarrow{}$\ignorespaces\fi}
\def\iff{\ifmmode\Longleftrightarrow \else
     \unskip${}\Longleftrightarrow{}$\ignorespaces\fi}

\let\:=\colon
\newtheorem{Theorem}{Theorem}[section]
\newtheorem{Lemma}[Theorem]{Lemma}
\newtheorem{Corollary}[Theorem]{Corollary}
\newtheorem{Proposition}[Theorem]{Proposition}
\newtheorem{Remark}[Theorem]{Remark}

\newtheorem{Example}[Theorem]{Example}

\newtheorem{Definition}[Theorem]{Definition}

\let\epsilon\varepsilon
\let\phi=\varphi
\let\kappa=\varkappa
\def\qed{\ifhmode\textqed\fi
   \ifmmode\ifinner\quad\qedsymbol\else\dispqed\fi\fi}
\def\textqed{\unskip\nobreak\penalty50
    \hskip2em\hbox{}\nobreak\hfil\qedsymbol
    \parfillskip=0pt \finalhyphendemerits=0}
\def\dispqed{\rlap{\qquad\qedsymbol}}

\opn\dis{dis}
\def\pnt{{\raise0.5mm\hbox{\large\bf.}}}

\begin{document}

\title{Prime filtrations of monomial ideals and polarizations}
\author{Ali Soleyman Jahan}
\address{Ali Soleyman Jahan, Fachbereich Mathematik und
Informatik, Universit\"at Duisburg-Essen, 45117 Essen, Germany}
\email{ali.soleyman-jahan@stud.uni-duisburg- essen.de}

\date{}
\begin{abstract} We show that all monomial ideals in the
polynomial ring in at most 3 variables are pretty clean and  that
an arbitrary  monomial ideal $I$ is pretty clean if and only if
its polarization $I^p$ is clean. This yields  a new
characterization of pretty clean monomial ideals in terms of the
arithmetic degree, and it also implies that a multicomplex is
shellable if and only the  simplicial complex corresponding to its
polarization is (non-pure) shellable. We also discuss Stanley
decompositions in relation to prime filtrations.
\end{abstract}
\maketitle

\section*{Introduction}
Let $R$ be a Noetherian ring, and $M$ a finitely generated
$R$-module. A basic fact in commutative algebra \cite[Theorem
6.4]{M} says that there exists a finite filtration
\[
\mathcal{F}:\; 0=M_0 \subset M_1 \subset \ldots \subset M_{r}=M
\]
with cyclic quotients $M_{i}/M_{i-1}\iso R/P_{i}$ and $P_{i} \in
\Supp(M)$. We call any such filtration of $M$ a prime filtration.
The set of prime ideals ${P_{1},\ldots,P_{r}}$ which define  the
cyclic quotients of $\mathcal{F}$  will be denoted by
$\Supp(\mathcal{F})$. Another basic fact \cite[Theorem 6.5]{M}
implies that $\Ass(M ) \subset \Supp(\mathcal{F})\subset\Supp(M)$.
Let $\Min (M )$ denote the set of minimal prime ideals in
$\Supp(M)$. Dress \cite{D} calls a prime filtration $\mathcal{F}$
of $M$ {\em clean} if $\Supp(\mathcal{F})=\Min(M)$. The $R$-module
$M$ is called {\em clean} if it admits a clean filtration. Herzog
and Popescu \cite{HP} introduced the concept of {\em pretty clean
modules}.

 A prime filtration
\[
\mathcal{F }:\;  0=M_0 \subset M_1 \subset \ldots \subset M_{r}=M
\]
of $M$ with $M_{i}/ M_{i-1} \iso R/P_{i}$ is called {\em pretty
clean},  if for all $i<j$ for which $P_{i} \subseteq P_{j}$ it
follows that $P_{i}=P_{j}$.

 In  other words, a proper inclusion
$P_{i} \subset P_{j}$ is only possible if $i>j$. The module $M$ is
called {\em pretty clean}, if it has a pretty clean filtration. We
say an ideal $I \subset R$ is pretty clean if $R/I$ is pretty
clean.

A prime filtration which is pretty clean has the nice property
that $\Supp(\mathcal{F})=\Ass(M)$, see \cite[Corollary 3.6]{HP}.
It is still an open problem to characterize the modules which have
a prime filtration $\mathcal F$  with
$\Supp(\mathcal{F})=\Ass(M)$. In Section 4 we give an example of a
module which is not pretty clean but nevertheless has a prime
filtration whose support coincides with the the set of associated
prime ideals of $M$.

 Dress showed in his paper \cite{D} that a simplicial complex is shellable
if and only if its Stanley-Reisner ideal is clean, and Herzog and
Popescu generalized this result by showing that the multicomplex
associated with a monomial ideal $I$ is shellable if and only if
$I$  is pretty clean. As a main result of this paper we relate
these two results by showing in Theorem \ref{Ghomri} that a
monomial ideal is pretty clean if and only if its polarization
$I^p$ is clean. As a consequence of this result we are able to
give the following characterization (Theorem \ref{adeg}) of pretty
clean monomial ideals: for any monomial ideal $I$ the length of
any prime filtration is bounded below by the arithmetic degree of
$I$, and equality holds if and only if $I$ is pretty clean.

In the first section of this paper we show that all monomial
ideals in $K[x_1,\ldots, x_n]$ of height $\geq n-1$ are pretty
clean and use this fact to show that any monomial ideal in the
polynomial ring in three variables is pretty clean; see
Proposition \ref{n-1} and Theorem \ref{xyz}. However for all
$n\geq 4$ there exists a monomial ideal of height $n-2$ which is
not pretty clean, see Example \ref{khob}.

In Section 2 we discuss the Stanley conjecture concerning Stanley
decompositions. In \cite[Theorem 6.5]{HP} it was shown that the
Stanley conjecture holds for any pretty clean monomial ideal.
Therefore using the results of Section 1 we recover the result of
Apel \cite[Theorem 5.1]{A} that the Stanley conjecture holds for
any monomial ideal in the polynomial ring in three variables.
Similarly we conclude that the Stanley conjecture holds for any
monomial ideal of codimension 1.

We also notice (Proposition \ref{imp}) that for a monomial ideal,
instead of requiring that $I$ is pretty clean, it suffice to
require that there exists a prime filtration $\mathcal {F}$ with
$\Ass(S/I)=\Supp(\mathcal{ F})$ in order to conclude that the
Stanley conjecture holds for $S/I$.

Unfortunately it is not true that each Stanley decomposition
corresponds to a prime filtration as shown by an example of
MacLagan and Smith \cite[Example 3.8]{MS}. However we characterize
in Proposition \ref{ma} those Stanley decomposition of $S/I$ that
correspond to prime filtrations. Using this characterization we
show in Corollary \ref{stp} that in the polynomial ring in two
variables Stanley decompositions and prime filtrations are in
bijective correspondence.

In Section 3 we prove the above mentioned result concerning
polarizations.  One important step in the proof (see Proposition
\ref{martaba}) is to show that there is a bijection between the
facets of the multicomplex defined by the monomial ideal $I$ and
the facets of the simplicial complex defined by the polarization
 of   this shows that $I$ and  $I^p$ have the same arithmetic degree.

 The final  Section  is devoted to prove the new characterization
of pretty clean monomial ideals in terms of the arithmetic degree.

\medskip
I would like to thank Professor Jürgen Herzog for many helpful
comments and discussions.

\section{Pretty clean monomial ideals and multicomplexes}

We denote by $S=K[x_1,\ldots, x_n]$ the polynomial ring in $n$
variables over a field $K$. Let $I\subset S$ be a monomial ideal.
In this paper a prime filtration of $I$ is always assumed to be a
monomial prime filtration. This means a prime filtration
\[
\mathcal{F}:\; I=I_0\subset I_1\subset \ldots\subset I_r=S
\]
with $I_j/I_{j-1}\iso S/P_j$, for $j=1,\ldots,r$ such that all
$I_j$ are monomial ideals.

Recall that the prime filtration $\mathcal{F}$ is called {\em
 pretty clean},  if for all $i<j$ which $P_i\subseteq P_j$  it
follows that $P_i=P_j$. The monomial ideal $I$ is called {\em
pretty clean}, if it has a pretty clean filtration.

In this section we will show that monomial ideals in at most three
variables are pretty clean.

Let $I \subset S $ be a monomial ideal. The saturation $\tilde{I}$
of $I$ is defined to be
\[
\tilde{I}=I:\mm^\infty=\Union_k (I:\mm^k),
\]
where $\mm=(x_1,\ldots,x_n)$  is the graded  maximal ideal of $S$.

We first note the following

\begin{Lemma}
\label{i}
 Let $I \subset S $ be a monomial ideal of $S$. Then $
I$ is pretty clean if and only if $\tilde{I}$ is pretty clean.
\end{Lemma}
\begin{proof} The $K$-vectorspace  $\tilde{I}/I$ has a finite dimension,
and we can choose monomials $u_1,\ldots, u_t\in \tilde{I}$ whose
residue classes modulo $I$ form a $K$-basis of $\tilde{I}/I$.
Moreover the basis can be chosen such that for all $j=1,\ldots, t$
one has $I_j/I_{j-1}\iso S/\mm$ where $I_0=I$ and
$I_j=(I_{j-1},u_j)$, and where $\mm=(x_1,\ldots,x_n)$  is the
graded  maximal ideal of $S$. Indeed, we have $\tilde{I}=I:\mm^k$
for some $k$. For each $i\in[k]$, where $[k]=\{1,\ldots,k\}$,  the
$K$-vectorspace $(I:\mm^i)/(I:\mm^{i-1})$ has  finite dimension.
If
\[
\dim_K (I:\mm^i/I:\mm^{i-1})=r_i,
\]
then we can choose monomials $u_{i1},\ldots,u_{ir_i} \in I:\mm^i$
whose residue classes modulo $I:\mm^{i-1}$ form  a basis for this
$K$-vectorspace. Composing these bases we obtain the required
basis for $\tilde{I}/I$.

So we have
\[
{\mathcal F}_1 :\; I=I_0 \subset I_1 \subset \dots \subset
I_{t}=\tilde{I}
\]
with $I_{i}/I_{i-1}\iso S/\mm$, for all $i=1,\ldots, t$. Now if
$\tilde{I}$ is a pretty clean and ${\mathcal G}$ is pretty the
clean filtration of $\tilde{I}$, then the prime filtration
$\mathcal F$ which is obtained by composing  ${\mathcal F}_1$ and
$\mathcal G$ yields a pretty clean filtration of  $I$.

For the converse, let $I=I_0 \subset I_1 \subset \dots \subset
I_r=S$ be pretty clean filtration of $I$. We will show that
$\tilde{I}$ is pretty clean by induction on $\dim_K
\tilde{I}/I=t$.  If  $t=0$ the assertion is  trivially true.
Assume now that  $t>0$. It is clear that  $I_1$ is also pretty
clean and that $I_1/I\iso S/\mm$, since $I\neq \tilde{I}$. It
follows that  $\tilde{I_1}=\tilde{I}$ and that $\dim_K
\tilde{I_1}/I_1 =t-1$.  So by the induction hypothesis
$\tilde{I}=\tilde{I_1}$ is pretty clean.
\end{proof}

\begin{Corollary}
\label{n}
 Let $S=K[x_1, \ldots,x_{n}]$ be the polynomial ring in $n$
variables. Then any monomial ideal of height $n$ is pretty clean.
\end{Corollary}

Our next goal is to show that even  the monomial ideals in
$S=K[x_1, \ldots,x_{n}]$ of height $\geq n-1$ are pretty clean. To
this end we have to recall the concept of multicomplexes and
shellings.

Stanley \cite{Stan1} calls a subset $\Gamma \subseteq \NN^n$ a
{\em multicomplex} if for all $a \in \Gamma$ and for all $b \leq
a$ one has $b\in\Gamma$. Herzog and Popescu \cite{HP}  give the
following modification of Stanley's definition of multicomplex
which will be used in this paper. Before we give this definition
we introduce some notation. We  set $\NN_\infty=\NN\union
\{\infty\}$. Let $\Gamma$ be a subset of $\NN^n_\infty$. An
element $m\in\Gamma$ is called maximal if there is no $a \in
\Gamma$ with $a>m$. We denote by ${M}(\Gamma)$ the set of maximal
elements of $\Gamma$. If $a\in\Gamma$, we call
\[
\infpt(a)=\{i:\;a(i)=\infty\},
\]
the {\em infinite part} of $a$.

\begin{Definition}{\em A subset $ \Gamma \subset \NN_\infty^{n}$   is
called a {\em multicomplex } if
\begin{enumerate}
\item[(i)] for all $a \in  \Gamma $ and  for all $ b \leq a $ it follows
that $b \in \Gamma $,
\item[(ii)] for all $a\in\Gamma$ there exists an element $m\in M(\Gamma)$ such that $a\leq m$ .
\end{enumerate}}
\end {Definition}

The elements of a multicomplex are called {\em faces}. An element
$a\in\Gamma$ is called a {\em facet} of $\Gamma$ if for all $m\in
{M}(\Gamma)$ with $a\leq m$ one has $\infpt(a)=\infpt(m)$. The set
of all facets of $\Gamma$ will be denoted by ${F}(\Gamma)$. The
facets in ${M}(\Gamma)$ are called {\em maximal facets}. It is
clear that ${M}(\Gamma)\subset {F}(\Gamma)$. We recall that for
each multicomplex $\Gamma$ the set of facets of $\Gamma$
 is a finite set, see \cite[Lemma 9.6]{HP}.

 Let $\Gamma$ be a multicomplex, and let $I(\Gamma)$ be the $K$-vectorspace in
$S=K[x_1,\ldots,x_n]$ spanned by all monomials $x^a$ such that
$a\not\in\Gamma$. Note that $I(\Gamma)$ is a monomial ideal, and
called the monomial ideal associated to $\Gamma$.

Conversely  let $I\subset S$ be any monomial ideal, then there
exists a unique multicomplex $\Gamma(I)$ with $I(\Gamma(I))=I$.
Indeed, let $A=\{a\in\NN^n:x^a\not\in I\}$; then
$\Gamma(I)=\Gamma(A)$ is called the multicomplex associated to $I$
, where $\Gamma(A)=\{b\in\NN_\infty^n :\; b\leq a \quad\text{for
some}\quad a\in A \}$ .

A subset $S\subset \NN_\infty^n$ is called a {\em Stanley} set if
there exists $a\in\NN^n$ and $m\in\NN_\infty^n$ with
$m(i)\in\{0,\infty\}$ such that $S=a+{S}^*$, where
${S}^*=\Gamma(m)$.

In \cite{HP} the concept of {\em shelling} of multicomplexes was
introduced as in the following by Herzog and Popescu.

\begin{Definition} \em{A multicomplex $\Gamma$ is {\em shellable} if
the facets of $\Gamma$ can be ordered $a_1,\ldots,a_r$ such that
\begin{enumerate}
\item[(i)] $S_i=\Gamma(a_i)\backslash \Gamma(a_1,\ldots,a_{i-1})$
is a Stanley set for all $i=1,\ldots,r$, and
\item[(ii)] whenever
${S_i}^* \subseteq {S_j}^*$, then ${S_i}^*={S_j}^*$ or $i>j$.
\end{enumerate}}
\end{Definition}
Any order of the facets satisfying ($i$) and ($ii$) is called a
{\em shelling} of $\Gamma$.

 In \cite[Theorem 10.5]{HP} the following has been
proved.

\begin{Theorem} The multicomplex $\Gamma$ is shellable if and only
if $S/I(\Gamma)$ is a  pretty clean $S$-module.
\end{Theorem}

\begin{Remark}\label{ah}{\em Let $\Gamma\subset \NN^n_\infty$ be a shellable multicomplex with
shelling $a_1,\ldots,a_r$, then $a_1(i)\in\{0,\infty\}$ and
therefore $a_1$ is one  of the minimal elements in $F(\Gamma)$
with respect to its partially order. Indeed,  since
$a_1,\ldots,a_r$ is a shelling, it follows that $S_1=\Gamma(a_1)$
is a Stanley set and therefore there exists a vector $b\in\NN^n$
and a vector $m\in\{0,\infty\}^n$ such that
 \begin{eqnarray*}
 \label{first}
 \Gamma(a_1)=b+\Gamma(m).
\end{eqnarray*}
It is clear that $\infpt(a_1)=\infpt(m)$. If $\infpt(m)=[n]$, then
there is nothing to show. Suppose now that $\infpt(m)\neq [n]$,
and choose $i\in [n]\setminus \infpt(m)$. If $a_1(i)\neq 0$ there
exists $c\in \Gamma(a_1)$ with $c(i)<a_1(i)$. Since $c$ and
$a_1\in b+\Gamma(m)=\Gamma(a_1)$, and since $m(i)=0$, it follows
that $c(i)=b(i)=a_1(i)$, a contradiction.

Furthermore,  if $\Gamma$ has only one maximal facet,  then
$F(\Gamma)$ has only one minimal element, and any shelling of
$\Gamma$ must start with this minimal element and end by the
maximal one. In fact, suppose  $a_1$ and $a_2$ are  minimal
elements in $F(\Gamma)$. By  the first part of this remark it
follows that $a_1$ and $a_2$ are vectors in $\{0,\infty\}^n$.
Hence since $\infpt(a_1)=\infpt(a_2)$, we see that  $a_1=a_2$. Now
let $a_1,\ldots, a_r$ be any shelling of $\Gamma$. Then, by what
we have shown, it follows that $a_1$  is the unique minimal
element in $F(\Gamma)$.  Let $m$ be the maximal element of
$F(\Gamma)$. Suppose  $m=a_k$ for some  $k<r$, then
\[
S_{k+1}=\Gamma(a_{k+1})\setminus
\Gamma(a_1,\ldots,a_k)=\Gamma(a_{k+1})\setminus
\Gamma(m)=\emptyset,
\]
which is not a Stanley set, a contradiction.  Moreover in this
case for each $i$ there exists a $d_i\in\NN^n$ such that $
S_i=d_i+\Gamma(a_1)$.}
\end{Remark}
 Now we are ready to  show that in $S=K[x_1,\ldots,x_n]$, any
ideal of height $n-1$ is pretty clean.

\begin{Proposition}\label{n-1}  Let $I \subset S=K[x_1, \ldots,x_n]$
be any monomial ideal of height $\geq n-1$.  Then $I$ is pretty
clean.
\end{Proposition}

\begin{proof}
We may assume that $I$ is a monomial ideal of height $n-1$, and by
Lemma \ref{i} that $I$ is saturated, i.e, $I=\tilde{I}$. It
follows that $I=\bigcap I_j$, where $I_j=({x_1}^{c_{j1}}, \ldots,
{x}_{j-1}^{c_{jj-1}},
{x}_{j+1}^{c_{jj+1}},\ldots,{x_n}^{c_{jn}})$, and where $c_{jk}>0$
for $k\neq j$. We denote by $\Gamma$ and $\Gamma_j$ the
multicomplexes associated to $I$ and $I_j$, and  by ${F}$ and
${F}_j$ the sets of facets of $\Gamma$ and $\Gamma_j$,
respectively.  The sets  $F$ and $F_j$ are finite, see \cite[Lemma
9.6]{HP}. Suppose $|{F}|=t$ and $|{F}_j|=t_j$.  Since $I_j$ is
$P_j$-primary where $P_j= (x_1, \ldots, x_{j-1},
{x}_{j+1},\ldots,{x_n})$, it follows  from \cite[Proposition
5.1]{HP} that $I_j$ is pretty clean,  and hence $\Gamma_j$ is
shellable. Moreover $a\in \NN^n_\infty$ is a facet of $\Gamma_j$
if and only if $a(j)=\infty$ and $a(k)<c_{jk}$ for $k\neq j$. Let
$a_{j1},\ldots,a_{jt_j}$ be a shelling of $\Gamma_j$.

 For showing $I$ is pretty clean it is enough to show that $\Gamma$ is
shellable.

By \cite[Lemma 9.9 (b)]{HP} we have
$\Gamma=\bigcup_{j=1}^n\Gamma_j$. Also by \cite[Lemma 9.10]{HP},
each $F_j$ has only one maximal facet, say $m_j$, where
\[
m_j(k)=
\begin{cases}
\infty,\quad\text{if}\quad k=j,
\\c_{j_k}-1,\quad\text{otherwise}.
\end{cases}
\]
It follows that ${F}=\bigcup {F}_j$ and that the union is
disjoint, since $a\in F$  belongs to $F_j$ if and only if
$a(j)=\infty$ and $a(k)<\infty$ for $k\neq j$. In particular one
has $(\Union_{i=1}^{j-1} F_i)\cap F_j=\emptyset$ for
$j=2,\ldots,n$.

We claim that
$$a_{11},\ldots,a_{1t_1},a_{21},\ldots,a_{2t_2},\ldots,a_{n1},\ldots,a_{nt_n}$$
is a shelling for $\Gamma$. Indeed, for all $j$ and all $k$ with
$1<k\leq t_j$ we have
\[
S_{jk}=\Gamma(a_{jk})\backslash\Gamma(a_{11},\ldots,a_{jk-1})=\Gamma(a_{jk})
\backslash \Gamma(a_{j1},\ldots,a_{jk-1}),
\]
and if $k=1$, then
\[
S_{j1}=\Gamma(a_{j1})\backslash\Gamma(a_{11},\ldots,a_{j-1t_{j-1}})=\Gamma(a_{j1}).
\]
Since $a_{j1},\ldots,a_{jt_j}$ is a shelling of $\Gamma_j$, it
follows that  $S_{jk}$  is a Stanley set for all $j$ and all $k$.

Condition $(ii)$ in the  definition of shellability is obviously
satisfied. In fact, since $\Gamma_j$ is shellable and has only one
maximal facet, it follows by Remark \ref{ah} that  for all
$k=1,\ldots,t_j$, there exists some $d_{jk}\in\NN^n$ such that
$S_{jk}=d_{jk}+S_j^*$, where $S_j^*=\Gamma(a_{j1})$. Moreover if
 $j\neq t$ then $ a_{j1}$ and $a_{t1}$ are not comparable,
and  hence in this case  there is no inclusion among $S_j^*$  and
$S_t^*$.
\end{proof}

As a Consequence of  Proposition \ref{n-1} we have

\begin{Corollary} \label{xy} Any monomial ideal $I\subset S=K[x,y]$ is pretty clean.
\end{Corollary}

Next we will show that any monomial ideal in $S=K[x_1,x_2,x_3]$ is
also  pretty clean. First we need

\begin{Lemma}
\label{ui} Let $I\subset S=K[x_1,x_2,x_3]$ be a monomial ideal of
height $1$. Then   $I=uJ$, where $u$ is a monomial in $S$, and $J$
is a monomial ideal of height $\geq 2$. Moreover, $I$ is pretty
clean if and only if $J$ is pretty clean.
\end{Lemma}

\begin{proof} The first statement of the lemma is obvious.
Assume now that $J$ is pretty clean with pretty clean
filtration
\[
\mathcal{F}:\; J=J_0 \subset J_1 \subset \ldots \subset J_r=S
\]
such that $J_i/J_{i-1} \iso S/P_i$, where $P_i\in\Ass J$. Then
$\height P_i\geq 2$. It follows that
\[
\mathcal{F}_1:\; I=uJ \subset uJ_1 \subset \ldots \subset uJ_r=(u)
\]
is a prime filtration of $(u)/I$ with factors $uJ_i/uJ_{i-1} \iso
S/P_i$.

There exists a prime filtration
\[
\mathcal{F}_2:\; (u)=J_r \subset J_{r+1}\subset \ldots \subset
J_{r+t}=S
\]
of the principal monomial ideal $I_1=(u)$, where the $J_{r+k}$ are
again  principal monomial ideals with $J_{r+k}/J_{r+k-1} \iso
S/Q_k$ and where $Q_k\in \Ass (u)$ has height $1$ for all $k$. In
fact, if $u=u_0=\prod_{t=1}^k x_{i_t}^{a_t}$ and
$u_j=\prod_{r=j+1}^k x_{i_r}^{a_r}$ for $j=1,\ldots,k-1$, then the
prime filtration $\mathcal{F}_2$ is the following:
\[
\mathcal{F}_2 :\; J_r=(u) \subset (x_{i_1}^{a_1-1}u_1)\subset
(x_{i_1}^{a_1-2}u_1)\ldots \subset (u_1) \subset
(x_{i_2}^{a_2-1}u_2) \subset\ldots\subset (u_2)\subset \ldots
\subset (x_{i_k}) \subset S.
\]
Therefore this filtration  of $I_1=(u)$ is pretty clean. Now
composing the above  filtrations $\mathcal{F}_1 $ and
$\mathcal{F}_2$ we obtain  a pretty clean filtration of $I$.

The converse follows from Proposition \ref{n-1} , because
$\height(J)\geq 2$.
\end{proof}

 Combining Lemma \ref{ui} with Proposition \ref{n-1} we get
\begin{Theorem}\label{xyz} Any monomial ideal in a polynomial ring in at most three variables is pretty
clean.
\end{Theorem}

The following example shows that this theorem  can not be extended
to polynomial rings in more than three variables, and it also
shows that monomial ideals of height $<n-1$ may  not be pretty
clean.

\begin{Example}\label{khob}{\em Let $n=4$, and $\Gamma$ be the multicomplex
with facets $(\infty,\infty,0,0)$ and $(0,0,\infty,\infty)$. Then
$\Gamma$ is not shellable, and so the monomial ideal
$$I(\Gamma)=(x_1x_3,x_1x_4,x_2x_3,x_2x_4)\subset
K[x_1,x_2,x_3,x_4]$$  is not pretty clean.

More generaley, let $n>3$ and $a=(0,0,\infty,\ldots,\infty)$ and
$b=(\infty,\infty,0,\ldots,0)$ be two elements in
$\NN^n_{\infty}$. Then $\Gamma=\Gamma(a,b)$ is not a shellable
multicomplex, hence $I=(x_1,x_2)\sect(x_3,\ldots,x_n)$ is a square
free monomial ideal in $S=K[x_1,\ldots,x_n]$ which is not clean.}
\end{Example}

\section{Prime filtrations and Stanley decompositions}

Let $I\subset S=K[x_1,\ldots,x_n]$ be a monomial ideal, any
decomposition of $S/I$ as a direct sum of $K$-vectorspaces of the
form $uK[Z]$ where $u$ is a monomial in $K[X]$, and $Z \subset
X=\{x_1,\ldots,x_n\}$ is called a {\em Stanley decomposition}. In
this paper we will call $uK[Z]$ a Stanley space of dimension
$|Z|$, where $|Z|$ denotes the cardinality of $Z$. Stanley
decomposition have been studied in various combinatorial and
algebraic contexts, see \cite{A},\cite{HT}, and \cite{MS}.

Let $R$ be a finitely generated standard graded $K$-algebra where
$K$ is a field, and let $M$ be a finitely generated graded
$R$-module.  Then the Hilbert series of $M$ is defined to be
$\Hilb(M)=\sum_{i\in Z}\dim_K M_it^i$. It is known that if
$\dim(M)=d$, then there exists a $Q_M(t)\in Z[t,t^{-1}]$ such that
 \[
 \Hilb(M)=Q_M(t)/(1-t)^d
 \]
 and $Q_M(1)\neq 0$. The number $Q_M(1)$ is called the
 multiplicity of $M$, and is denoted by $\e(M)$.

Let $I\subset S$ be a monomial ideal. Then the number of Stanley
spaces of a given dimension in a Stanley decomposition may depend
on this particular decomposition. For example, if $I=(xy)\subset
K[x,y]$, then for all integers $k>0$ and $l>0$ we have the Stanley
decomposition
$$S/I=x^lK[x]\dirsum y^kK[y]\dirsum
(\Dirsum_{i=0}^{l-1}x^iK)\dirsum (\Dirsum_{j=1}^{k-1}y^jK),$$
 for $S/I$ with as many Stanley spaces of dimension 0 as we want,
however only 2 Stanley spaces of dimension 1 in any Stanley
decomposition.

This is a  general fact. Indeed, the number of Stanley spaces of
maximal dimension is independent of the special Stanley
decomposition. In fact, this number is equal to the multiplicity,
$e(S/I)$, of $S/I$.

 Let \[
  S/I=\Dirsum_{i=1}^r u_iK[Z_i]
  \]
  be an arbitrary Stanley decomposition of $S/I$, and $d=\max\{|Z_i|\:i=1,\ldots,r\}$. Then
  \[
  \Hilb(S/I)=\sum_{i=1}^r\Hilb(u_iK[Z_i])=\sum_{i=1}^r
  t^{\deg(u_i)}/(1-t)^{|Z_i|}=Q_{S/I}(t)/(1-t)^d.
  \]
  with  $Q_{S/I}(t)=\sum_{i=1}^r
  (1-t)^{d-|Z_i|}t^{\deg(u_i)}$.  It follows that  $e(S/I)=Q_{S/I}(1)$ is equal to the
  number of Stanley space of dimension $d$ in this Stanley
  decomposition of $S/I$.

We also note that for each monomial $u\in\tilde{I}\setminus I$ the
0-dimensional Stanley space $uK$ belongs to any Stanley
decomposition of $S/I$. In fact $u\mm^k\subset I$ for some $k$.
Now if $u$ belongs to some Stanley space $vK[Z]$ with $|Z|\geq 1$,
then $vK[Z]\cap I\neq\emptyset$, a contradiction.

\medskip
Stanley \cite{Stan} conjectured that there always exists a Stanley
decomposition $$S/I=\Dirsum_{i=1}^r u_iK[Z_i],$$ such that
$|Z_i|\geq\depth(S/I)$ for all $i$.

\medskip
Apel \cite{A} studied some cases in which Stanley's conjecture
holds. Also in \cite[Theorem 6.5]{HP} it has been proved that for
all pretty clean monomial ideals Stanley's conjecture holds.
Therefore combining  Theorem \ref{xyz} and Lemma \ref{n-1} with
\cite[Theorem 6.5]{HP} we get

\begin{Proposition}\label{apel}
{\em (a)} Let $I\subset S=K[x_1,\ldots,x_n]$ be a monomial ideal
of height $\geq n-1$. Then Stanley's conjecture holds for $S/I$.

{\em (b) (Apel, \cite[Theorem 5.1]{A})} Let $I$ be any monomial
ideal in the polynomial ring in at most three variables. Then
Stanley's conjecture holds for $S/I$.
\end{Proposition}

In the proof of  \cite[Theorem 6.5]{HP} it is used that Stanley
decompositions of $S/I$ arise from prime filtrations. In fact, if
$\mathcal F$  is a  prime filtration of $S/I$ with factors
$(S/P_i)(-a_i)$, $i=1,\ldots,r$. Then if we set $u_i=\prod_{j=1}^n
x_j^{a_i(j)}$ and $Z_i=\{x_j :\; x_j\not\in P_i\}$, then
\[
S/I=\Dirsum_{i=1}^r u_iK[Z_i].
\]

Recall that if $\mathcal{F}$ is a pretty clean filtration of
$S/I$, then $\Ass(S/I)=\Supp(\mathcal{F})$. The converse of this
statement is not always true, see Example \ref{gut}. As a
generalization of \cite[Theorem 6.5]{HP} we show

\begin{Proposition}\label{imp}Suppose $I\subset S$ is a monomial ideal, and
$\mathcal{F}$ is a  prime filtration of $S/I$ with
$\Supp(\mathcal{F})=\Ass(S/I)$. Then the Stanley decomposition of
$S/I$ which is obtained from this prime filtration satisfies the
condition of Stanley's conjecture.
\end{Proposition}

\begin{proof} The Stanley decomposition which is obtained from
$\mathcal{F}$ has the property that $|Z_i|=\dim S/P_i$. By
\cite[Proposition 1.2.13]{BH} we have $\depth(S/I)\leq\dim(S/P_i)$
for all $P_i\in\Ass(S/I)$, and hence the assertion follows.
\end{proof}

In all cases discussed above we found a Stanley decomposition
corresponding  to a prime filtration and  satisfying  the Stanley
conjecture. However we will  show that there exist  examples of
monomial ideals such that {\em all} Stanley decompositions arising
from a prime filtration may fail to satisfy the Stanley
conjecture.

 First we notice that
 \begin{Remark}\label{plast} {\em Let $I\subset S$ be a Cohen-Macaulay monomial ideal, and
 $$\mathcal{F}:\; I=I_0\subset I_1\subset\ldots\subset I_r=S$$
 be a prime filtration of $S/I$. We claim that if the Stanley
 decomposition of $S/I$ corresponding  to $\mathcal{F}$ satisfies the
 Stanley conjecture, then $\Ass(I)=\Supp(\mathcal{F})$. In particular $I$ is clean, since $\Min(I)=\Ass(I)$.

Indeed,
 since $I$ is Cohen-Macaulay we have
 $\depth(S/I)=\dim(S/I)=\dim(S/P)$ for all $P\in\Ass(I)$.  We recall that
  $I_i/I_{i-1}\iso S/P_i(-a_i)$ for suitable $a_i$ and that $P_i\in\Ass(I_{i-1})$ for
  $i=1,\ldots,r$. Let  $T_i=u_iK[Z_i]$ be the Stanley space
  corresponding to $S/P_i(-a_i)$ as explained as above. Then  $|Z_i|=\dim(S/P_i)$. Assume that  $P_i\not\in\Ass(I)$ for some
  $i>1$. Since $I\subset I_{i-1}\subset P_i$, there exists a
  $P_j\in\Ass(I)$ such that $P_j\subsetneq P_i$. It follows that
  $|Z_i|=\dim(S/P_i)<\dim(S/P_j)=\depth(S/I)$, a contradiction.}
  \end{Remark}

  \begin{Example}\label{last} {\em Let $K$ be a field and $$I=
  (abd, abf, ace, adc, aef, bde, bcf, bce, cdf,
def)\subset
  S=K[a,b,c,d,e,f].$$
 The ideal $I$ is the Stanley-Reisner ideal corresponding  to the  simplicial complex
  $\Delta$ which  is the triangulation of the real
  projective plane $\PP^2$, see \cite[Figure 5.8]{BH}.
    It is  known that $S/I$ is Cohen-Macaulay
   if and only if $\chara(K)\neq 2$.
    This implies  $S/I$ is not  clean, since otherwise $\Delta$ would be  shellable
    and $S/I$ would be  Cohen-Macaulay for any field $K$. Hence  by Remark
    \ref{plast}, if $\chara(K)\neq 2$,
    no Stanley decomposition of $S/I$ which  corresponds to a prime filtration
     of $S/I$  satisfies the Stanley conjecture.}
     \end{Example}

Unfortunately not all Stanley decompositions of $S/I$  correspond
to prime filtrations, even if $S/I$ is pretty clean. Such an
example is given by McLagan and Smith in \cite{MS}. Let
$I=(x_1x_2x_3)\subset K[x_1,x_2,x_3]$. Then
$$S/I=1\dirsum x_1K[x_1,x_2]\dirsum x_2K[x_2,x_3]\dirsum
x_3K[x_1,x_3]$$ is a Stanley decomposition of $S/I$ which does not
correspond to a prime filtration of $S/I$. On the other hand, by
Theorem \ref{xyz} we know that $S/I$ is a pretty clean.

Now we want to characterize those Stanley decompositions of $S/I$
which correspond to a prime filtration of $S/I$.

We fist notice

\begin{Lemma}\label{her} Let $I\subset S=K[x_1,\ldots,x_n]$ be a monomial
ideal, and $T= uK[Z]$ be a Stanley space in a Stanley
decomposition of $S/I$. The $K$-vectorspace $I_1=I\dirsum T$ is a
monomial ideal if and only if $I_1=(I,u)$. In this case, $I:u=P$,
where $P=(x_i\:x_i\not\in Z)$.
\end{Lemma}

\begin{proof} We have
$I\subset I_1$ and $u\in I_1$. Suppose now that $I_1$ is a
monomial ideal. Since $(I,u)$ is the smallest monomial ideal that
contains $I$ and $u$, it follows that $(I,u) \subset I_1$. On the
other hand,   $I_1=I+uK[Z]\subset I+uK[x_1,\ldots,x_n]=(I,u)$.
Hence $I_1=(I,u)$.

Since for each $x_i\not\in Z$ we have $x_iu\in I_1=I\dirsum T$ and
$x_iu\not\in uK[Z]=T$, it follows that $x_iu\in I$ and hence
$x_i\in I:u$. On the other hand, if $v\in K[Z]$ is a monomial,
then $vu\not\in I$, since $uK[Z]$ is a Stanley space of $S/I$.
Therefore  $I:u=P=(x_i\:x_i\not\in Z)$.
\end{proof}

 \begin{Corollary} The monomial ideal $I\subset S$ is a prime
 ideal if and only if there exists a Stanley decomposition of $S/I$
 consisting of only one Stanley space.
 \end{Corollary}

As a consequence of this Lemma we have

\begin{Proposition}\label{ma} Let $I\subset S$ be a monomial ideal, and
$S/I=\Dirsum_{i=1}^ru_iK[Z_i]$ be a Stanley decomposition of
$S/I$. The given Stanley decomposition corresponds to a prime
filtration of $S/I$ if and only if the Stanley spaces
$T_i=u_iK[Z_i]$ can be ordered $T_1,\ldots,T_r$, such that
$$I_k=I\dirsum T_1\dirsum\ldots\dirsum T_k$$
is a monomial ideal for $k=1,\ldots,r$.
\end{Proposition}

\begin{proof} We prove "if" by induction on $r$. If $r=0$  then the assertion is trivially true.
 Let $r\geq 1$. By assumption $I_1=I\dirsum T_1$ is a monomial ideal.
 Hence by Lemma \ref{her} we have $I_1=(I,u_1)$ and
 $I:u_1=P_1=(x_i\:x_i\not\in Z_1)$. We notice that in this case $I_1/I\iso
 S/P_1(-a_1)$ and $u_1=\prod_{j=1}^n x_j^{a_1(j)}$, and that
 $S/I_1=\Dirsum_{i=2}^rT_i$. Now by the induction hypothesis this
 Stanley decomposition of $S/I_1$ corresponds to a prime
 filtration, say $\mathcal{F}_1$
 \[
 \mathcal{F}_1:\; I_1\subset I_2\subset\ldots\subset I_r=S.
 \]
 Therefore the given Stanley decomposition of $S/I$ corresponds to
 the prime filtration
\[
 \mathcal{F}:\; I\subset I_1\subset I_2\subset\ldots\subset I_r=S.
 \]
The converse follows immediately if we order the Stanley spaces of
$S/I$ which are  obtained from a prime filtration according to the
order of the ideals in this filtration.
\end{proof}

 We conclude this section by showing
\begin{Corollary}
\label{stp} Let $I\subset S=K[x,y]$ be a monomial ideal. Then each
Stanley decomposition of $S/I$ corresponds to a prime filtration
of $S/I$.
\end{Corollary}

\begin{proof}  The $K$-vectorspace $\tilde{I}/I$ has finite dimension, say $m$. So
we can choose monomials $v_1,\ldots, v_m\in \tilde{I}$ whose
residue classes modulo $I$ form a $K$-basis for $\tilde{I}/I$. As
observed in the discussions before Proposition \ref{apel}, in any
Stanley decomposition of $S/I$ these monomials have to appear as
$0$-dimensional Stanley spaces. In the proof of Lemma \ref{i} we
showed that it is possible to order the monomials $v_1,\ldots,
v_m$ in such a way that
\[
I_i=I\dirsum v_1K\dirsum\ldots\dirsum v_iK=(I,v_1,\ldots,v_i)
\]
 is a monomial ideal for $i=1\ldots,m$. If we remove in the given
Stanley decomposition of $S/I$ the Stanley spaces $v_iK$,
$i=1,\ldots,m$, the remaining summands establish a Stanley
decomposition of $S/\tilde{I}$. Thus we may assume that $I$ is
saturated. Hence $I=(x^\alpha y^\beta)$.

Let $S/I=\Dirsum_{i=1}^r u_iK[Z_i]$ be a Stanley decomposition of
$S/I$. We will prove  by induction on $\alpha+\beta$ that the
given Stanley decomposition can be ordered such that
$I_k=I\dirsum(\Dirsum_{i=1}^ku_iK[Z_i])$ is a monomial ideal for
all $k$. If $\alpha+\beta=0$ the assertion is trivially true. Let
$\alpha+\beta>0$. The Stanley decomposition of $S/I$ contains at
least one summand of the form $x^{\alpha-1}y^{\gamma}K[y]$, where
$\gamma\geq\beta$, or $x^{\theta}y^{\beta-1}K[x]$, where
$\theta\geq\alpha$.

We may assume that $x^{\alpha-1}y^{\gamma}K[y]$ is one of the
summands.   Let $t=\gamma-\beta$, and set
$v_i=x^{\alpha-1}y^{\gamma-i+1}$ for $i=1,\ldots,t+1$. If we set
$T_1=v_1K[y]$, then $I_1=I\dirsum T_1=(I,v_1)$ is a monomial
ideal. If we remove the Stanley space $T_1$ from the given Stanley
decomposition of $S/I$, the remaining establish a Stanley
decomposition of $S/I_{1}$. Since $v_2,\ldots,v_{t+1}$ belong to
$\tilde{I}_1\setminus I_1$,  these monomials have to appear in any
Stanley decomposition of $S/I_1$ as $0$-dimensional Stanley
spaces. In particular these monomials appear as 0-dimensional
Stanley space, $T_2=v_2K,\ldots,T_{t+1}=v_{t+1}K$ in the given
Stanley decomposition of $S/I$. Now it is clear that
$I_i=I_{i-1}\dirsum T_i=(I_{i-1},v_i)$ is a monomial ideal for
$i=1,\ldots,t+1$, where $I_0=I$.

Removing the Stanley spaces $T_1,\ldots,T_{t+1}$ from the given
Stanley decomposition of $S/I$, the remaining summands establish a
Stanley decomposition of $S/I_{t+1}$. Since
$I_{t+1}=(x^{\alpha-1}y^{\beta})$ is a saturated ideal, the
assertion follows by the induction hypothesis applied to
$S/I_{t+1}$.
\end{proof}

\section{A characterization of pretty clean monomial ideals in terms of polarizations}

In this section we  consider polarizations of monomial ideals and
of prime filtrations. Let $S=K[x_1,\ldots,x_n]$ be the polynomial
ring in $n$ variables over the field $K$, and $u=\prod_{i=1}^n
x_i^{a_i}$ be a monomial in $S$. Then
$$u^p=\prod_{i=1}^n\prod_{j=1}^{a_i}x_{ij}\in K[x_{11},\ldots,x_{1a_1},\ldots,x_{n1},\ldots,x_{na_n}]$$ is
called the {\em polarization} of $u$.

Let  $I$ be a monomial ideal in $S$ with monomial generators
$u_1,\ldots, u_m$. Then $(u_1^p,\ldots,u^p_m)$ is called {\em a
polarization of $I$}. Note that if $v_1,\ldots, v_k$ is a another
set of monomial generators of $I$ and if $T$ is the polynomial
with sufficiently many variables $x_{ij}$ such that  all the
monomials $u_i^p$ and $v_j^p$ belong to $T$, then
\[
(u_1^p,\ldots,u^p_m)T=(v_1^p,\ldots, v_k^p)T.
\]
Therefore we denote any polarization of $I$  by $I^p$, since in a
common polynomial ring extension all polarizations are the same,
and we write $I^p=J^p$ if a polarization of $I$ and a polarization
of $J$ coincide in a common polynomial ring extension.

Now let $I=(u_1,\ldots,u_m)\subset S$  be  a monomial ideal, and
$u\in S$ a monomial. Furthermore let  $T$ be the polynomial ring
in variables $x_{ij}$ such that:
\begin{enumerate}
\item[(1)]for all $i\in[n]$ there exists $k_i\geq 1$ such that
$x_{i1},\ldots,x_{ik_i}$ are in $T$, \item[(2)] $I^p\subset T$,
and $u^p\in T$.
\end{enumerate}
We consider the $K$-algebra homomorphism
$$\pi :T \longrightarrow S,\quad x_{ij} \mapsto x_i.$$
Then $\pi$ is an epimorphism with $\pi(u^p)=u$ for all monomials
$u\in S$, and $u^p$ is the unique squarefree monomial in $T$ of
the form $\prod_{i=1}^n\prod_{j=1}^{t_i} x_{ij}$ with this
property. In particular, $\pi(I^p)=I$. We call $\pi$ the
specialization map attached with the polarization.

\begin{Remark}\label{si}{\em Let $I=(u_1,\ldots,u_m)\subset S$ be a monomial ideal, and $u\in S$
a monomial. Then
\begin{enumerate}
\item[(a)] $I:u=(u_i/\gcd(u_i,u))_{i=1}^m$, and it is again a
monomial ideal in $S$.

\item[(b)] $I:u$ is a prime ideal if and only if for each $i\in[m]$,
  there exists a $j\in[m]$ such
 that $u_j/\gcd(u_j,u)$ is a monomial of  degree one, and
 $u_j/\gcd(u_j,u)$ divides $u_i/\gcd(u_i,u)$.

\item[(c)] Let $u=\prod_{i=1}^n x_i^{a_i}$ and $u_j=\prod_{i=1}^n
x_i^{b_i}$. If $u_j/\gcd(u_j,u)=x_i$, then $b_i=a_i+1$ and
$b_t\leq a_t$ for all $t\neq i$. Therefore $u_j/\gcd(u_j,u)=x_i$
if and only if $u_j^p/\gcd(u_j^p,u^p)=x_{ib_i}$.
\end{enumerate}
}
\end{Remark}

\begin{Lemma}
\label{ij} Let $I=(u_1,\ldots,u_m)\subset S$ be a monomial ideal
and $u\in S$ a monomial. If $I^p:u^p$ is a prime ideal,  then
$I^p:u^p=(x_{i_1j_1},\ldots,x_{i_kj_k})$ with $i_r\neq i_s$ for
$r\neq s$.
\end{Lemma}

\begin{proof} Since $I^p:u^p$ is a monomial prime ideal in polynomial ring $T$ it must be generated by variables.
If $x_{ij}$ and $x_{ik}$ are two generators of $I^p:u^p$, then
there exist $ r_j\in[m]$, and $ r_k\in[m]$ such that
$x_{ij}=u_{r_j}^p/\gcd(u_{r_j}^p,u^p)$ and
$x_{ik}=u_{r_k}^p/\gcd(u_{r_k}^p,u^p)$. It follows from Remark
\ref{si}(c) that $j-1=k-1$ is equal to the  exponent of $x_i$ in
$u$. Hence $x_{ij}=x_{ik}$.
\end{proof}

We also need to show
\begin{Lemma}
\label{ali} Let $I=(u_1,\ldots,u_m)\subset S$ be a monomial ideal,
and $u\in S$ a monomial in $S$. Then $I:u$ is a prime ideal if and
only if $I^p:u^p$ is a prime ideal.  In this case
$I:u=\pi(I^p:u^p)$.
\end{Lemma}

\begin{proof} Let $I:u$ be a prime
ideal. We may assume that $I:u=(x_1,\ldots,x_k)$ for some
$k\in[n]$. Therefore for each $i\in[k]$ there exists some
$u_{j_i}$, with $j_i\in[m]$, such that
$x_i=u_{j_i}/\gcd(u_{j_i},u)$ and for each $t\in[m]$, there exists
$i_t\in[k]$, such that $x_{i_t}$ divides $(u_t/\gcd(u_t,u))$.
Therefore by Remark \ref{si}(c) we have
$u_{j_i}^p/\gcd(u_{j_i}^p,u^p)=x_{it_i}$, where $t_i$ is the
exponent of $x_i$ in $u_{j_i}$ and $t_i-1$ is the exponent of
$x_i$ in $u$.

Also for each $s\in[m]$,  the monomial $u_s^p/\gcd(u_s^p,u^p)$ is
divided by one of thess $x_{it_i}$, where $ i\in[k]$. Indeed,
since $I:u$ is a prime ideal there exists some $i\in[k]$ such that
$x_i$ divides $(u_s/\gcd(u_s,u))$, where
$x_i=u_{j_i}/\gcd(u_{j_i},u)$. Let $t_i-1$ be the exponent of
$x_i$ in $u$. Then it follows  that the exponent of $x_i$ in $u_s$
is $> t_i-1$. Hence $x_{it_i}$ divides $u_s^p/\gcd(u_s^p,u^p)$,
and $I^p:u^p=(x_{1t_1},\ldots,x_{kt_k})$.

For the converse, let $I^p:u^p$ be a prime ideal. By Lemma
\ref{ij} we may assume that $I^p:u^p=(x_{1t_1},\ldots,x_{kt_k})$.
This means that for each $i\in [k]$ there is a monomial $u_{j_i}$
with $ j_i\in[m]$ such that
$x_{it_i}=u_{j_i}^p/\gcd(u_{j_i}^p,u^p)$ and for each $s\in [m]$,
the squarefree monomial $u_s^p/\gcd(u_s^p,u^p)$ is divided by one
of these $x_{it_i}$. Therefore by Remark \ref{si}(c) we have
$x_i=u_{j_i}/\gcd(u_{j_i},u)$ for $i\in[k]$, and for each $s\in
[m]$,  one of these variables divides $u_s/\gcd(u_s,u)$. Hence
$I:u=(x_1.\ldots,x_k)$.
\end{proof}

Let $I\subset S$ be a monomial ideal and
\[
\mathcal{F}:\; I=I_0\subset I_1\subset\ldots\subset I_r=S
\]
a  filtration of $S/I$.  We call $r$ the {\em length of
filtration} $\mathcal{F}$ and denote it by $\ell(\mathcal{F})$.

\medskip
Assume now that for all $j$ we have $I_{j+1}=(I_j,u_j)$ where
$u_j\in S$ is a monomial. We will define  {\em the polarization
$\mathcal{F}^p$} of $\mathcal{F}$ inductively as follow: set
$J_0=I^p$; assuming that $J_i$ is already  defined, we set
$J_{i+1}=(J_i,u_i^p)$. So $J_i=(I^p,u_1^p,\cdots, u_i^p)$,  and

\[
\mathcal{F}^p :\; I^p=J_0\subset J_1\subset\ldots\subset J_r=T
\]
is a filtration of $T/I^p$.

 We have the following
\begin{Proposition}\label{ffp} Suppose $I\subset S$ is a monomial ideal, and
\[
\mathcal{F}:\; I=I_0\subset I_1\subset\ldots\subset I_r=S
\]
a filtration of $S/I$ as  above. Then $\mathcal{F}$ is a prime
filtration of $S/I$ if and only if $\mathcal{F}^p$ is a prime
filtration of $T/I^p$.
\end{Proposition}
\begin{proof} Let
\[
\mathcal{F}:\; I=I_0\subset I_1\subset\ldots\subset I_r=S
\]
be a prime filtration of $S/I$. We use induction on
$r=\ell(\mathcal{F})$ the length of prime filtration. If $r=1$,
then $I$ is a monomial prime ideal and $I^p=I$.

Let $r> 1$. Then $\mathcal{F}_1:I_1\subset\ldots\subset I_r=S$ is
a prime filtration of $S/I_1$, and $\ell(\mathcal{F}_1)=r-1$.  By
our induction hypothesis, $\mathcal{F}_1^p$ is a prime filtration
of $I_1^p=(I^p,u_1^p)$. Since $I_1/I\iso I_1:u_1$ is a prime
ideal, it follows from Lemma \ref{ali} that  $J_0/J_1\iso
I_1^p:u_1^p$ is
 a prime ideal too. Hence $\mathcal{F}^p$ is a prime filtration of
$T/I^p$.

The other direction of the statement is proved similarly.
\end{proof}

Let $S=K[x_1,\ldots,x_n]$ a the polynomial ring, and $u,v\in S$ be
monomials. We notice that  $$\lcm(u,v)^p=\lcm(u^p,v^p).$$

Therefore we have
\begin{Lemma}
\label{ijp} Let $I,J$ be two monomial ideals in $S$. Then $(I\cap
J)^p=I^p\cap J^p$.
\end{Lemma}
\begin{proof} Let $I=(u_1,\ldots,u_m)$ and $J=(v_1,\ldots,v_t)$.
Then $I\cap J=(\lcm(u_i,v_j))$, where $1\leq i\leq m$ and $1\leq
j\leq t$.  Therefore $(I\cap J)^p=
(\lcm(u_i,v_j)^p)=(\lcm(u_i^p,v_j^p))= I^p\cap J^p$.
\end{proof}
We recall that a monomial ideal  $I\subset S$ is an irreducible
monomial ideal if and only if  there exists a subset $A\subset
[n]$ and for each $i\in A$ an integer $a_i>0$ such that
$I=(x_i^{a_i}\: i\in A)$, see \cite[Theorem 5.1.16]{Vil}. It is
known that for each minimal ideal $I$ there exists a decomposition
$I=\Sect_{i=1}^r J_i$ such that $J_i$ are irreducible monomial
ideals.

\begin{Corollary}\label{minI^p} Suppose $J_1,\ldots ,J_r$ are monomial ideals in the
polynomial ring $S$, and $I=\Sect_{i=1}^r J_i$. Then $I^p=
\Sect_{i=1}^r J_i^p$. In particular the minimal prime ideals of
$I^p$ are of the form $(x_{i_1t_1},\ldots,x_{i_kt_k})$, with
$i_r\neq i_s$ for $r\neq s$ .
\end{Corollary}

\medskip
Next we show that if $I\subset S$ is a monomial ideal and $I^p$
the polarization of $I$, then $|F(\Gamma(I))|=|F(\Gamma(I^p))|$.
First we notice the following:

\begin{Lemma}\label{H} Let $I\subset S$ be an irreducible monomial ideal
and $I^p$ the polarization of $I$. Furthermore, let $F$ and $F^p$
be the sets of facets of $\Gamma(I)$ and $\Gamma(I^p)$,
respectively. Then there exists a bijection between $F$ and $F^p$.
\end{Lemma}
\begin{proof} By \cite[Theorem 5.1.16]{Vil} there exists a subset
$A\subset [n]$ and for each $i\in A$ an integer $a_i>0$ such that
$I=(x_i^{a_i}\: i\in A)$. We may assume $A=[k]$ for some $k\leq
n$. In this case $\Gamma(I)=\Gamma(m)$, where
\[
m(i)=
\begin{cases}
a_i-1,\quad\text{if}\quad i\in[k],
\\ \infty,\quad\text{otherwise},
\end{cases}
\]
and $a\in F$ if and only if $a\leq m$ and $a(i)=\infty$ for $i>k$.
We have
$$I^p=(\prod_{j=1}^{a_1}x_{1j},\prod_{j=1}^{a_2}x_{2j},\ldots,\prod_{j=1}^{a_k}x_{kj}),$$
and we know that the facets in $F^p$ correspond to the minimal
prime ideals of $I^p$. Indeed, if $a\in F^p$ is a facet of
$\Gamma^p$, then $P_a=(x_i\:a(i)=0)$ is a minimal prime ideal of
$I^p$. Each minimal prime ideal of $I^p$ is of the form
$(x_{1t_1},\ldots,x_{kt_k})$, with $t_i\leq a_i$.

\medskip
Now we define
\[
\theta:F \to F^p, \; a\mapsto \bar{a}
\]
as follows: if $k<i\leq n$, then $\bar{a}(ij)=\infty$ for all $j$,
and if $i\in[k]$ we have $a(i)=t_i<a_i$, and we set
\[
\bar{a}(ij)=
\begin{cases}
0,\quad\text{if}\quad j=t_i+1,
\\ \infty,\quad\text{otherwise}.
\end{cases}
\]

Obviously $\bar{a}\in F^p$, since
$P_{\bar{a}}=(x_{1t_1+1},\ldots,x_{kt_k+1})$ is a minimal prime
ideal of $I^p$, and it is also clear  that $\theta$ is an
injective map.

 Let $\bar{a}\in F^p$. Then $\bar{a}$ corresponds to the minimal prime ideal
  $P_{\bar{a}}=(x_{1t_1},\ldots,x_{kt_k})$, where $t_i\leq a_i$. Therefore if $k<i\leq
  n$, we have  $\bar{a}(ij)=\infty$ for all $j$, and if $i\in[k]$, then
  \[
  \bar{a}(ij)=
  \begin{cases}
  0,\quad\text{if}\quad j=t_i,
  \\ \infty,\quad\text{otherwise}.
  \end{cases}
  \]
  Let $a\in\NN^n_{\infty}$  be the following:
  \[
  a(i)=
  \begin{cases}
  t_i-1,\quad\text{if}\quad i\in[k],
  \\ \infty,\quad\text{otherwise},
  \end{cases}
  \]
 then $a$ is a facet in $F$, since $a\leq m$ and $\infpt(a)=n-k=\infpt(m)$,  and moreover $\theta(a)=\bar{a}$.
\end{proof}

\medskip
Now let $I=(u_1,\ldots,u_m)\subset S$ be a monomial ideal and let
$D\subset [n]$ be the set of elements $i\in [n]$ such that $x_i$
divides $u_j$ for at least one $j=1,\ldots, m$.  Then we set
$$r_i=\max\{t\: \text{ $x_i^t$ divides $u_j$  at least for one  $j\in[m]$}\}$$
if $i\in D$ and $r_i=1$, otherwise. Moreover we set  $r=
\sum_{i=1}^nr_i$.

Note that $I$ has a decomposition $I=\Sect_{i=1}^t J_i$ where the
ideals $J_i$ are irreducible monomial ideals. In other words, each
$J_i$ is generated by pure powers of some of the variables. Then
$I^p= \Sect_{i=1}^tJ_i^p$ is an ideal in the polynomial ring
$$T=K[x_{11}\cdots, x_{1r_1}, x_{21}\cdots,\cdots, x_{n1},\cdots, x_{nr_n}]$$
in $r$ variables.

\medskip
We denote by $\Gamma$, $\Gamma^p$, $\Gamma_i$ and $\Gamma_i^p$ the
multicomplexes associated to $I$, $I^p$, $J_i$  and $J_i^p$,
respectively,
 and  by $F$, $F^p$, $F_i$ and $ F_i^p$ the sets of
facets of $\Gamma$, $\Gamma^p$, $\Gamma_i$ and $\Gamma_i^p$,
respectively.

It is clear that $F\subset \Union_{i=1}^t F_i$ since
$\Gamma=\Union _{i=1}^t\Gamma_i$, and also that
$F^p\subset\Union_{i=1}^t F_i^p$. Each $\Gamma_i$ has only one
maximal facet, say $m_i$, and $m_i(k)\leq r_k-1$ if
$m_i(k)\neq\infty$.

 let $A\subset \NN^n_{\infty}$ be the following
set:
$$A=\{a\in\NN^n_{\infty}\: a(i)<r_i \quad\text{if}\quad
a(i)\neq\infty\}.$$ We define the map
\[
\beta:A\to \{0,\infty\}^r,\;a \mapsto \bar{a}
\]
 as follows: if $a(i)=\infty$, then
 $\bar{a}(ij)=\infty$ for all $j$,
  and if $a(i)=e$ where $e\leq r_i-1$, then
  \[
  \bar{a}(ij)=
  \begin{cases}
  0,\quad\text{if}\quad j=e+1,\\ \infty,\quad\text{otherwise.}
  \end{cases}
  \]

 \begin{Proposition} \label{martaba} With the above assumptions and notation the restriction
 of the
  map $\beta$ to $F$ is a bijection from
 ${F}$ to ${F}^p$.

 \end{Proposition}

 \begin{proof} First of all we want to show that
 $\bar{a}\in {F}^p$. Indeed, $a\in {F}\subset \Union_{i=1}^t F_i$. Therefore there
exists an integer $j\in[n]$ such that $a\in F_j$, and since the
restriction of $\beta$ to $F_j$ is the map $\theta$ defined in
Lemma \ref{H}, it follows that  $\bar{a}\in F_j^p$. Therefore
there exists a subset $\{j_1,\ldots,j_s\}\subset [n]$ and positive
integers $t_k$ with $t_k\leq r_{j_k}$ for $k=1,\ldots,s$ such that
$P_{\bar{a}}=(x_{j_1t_1},\ldots,x_{j_st_s})$. It is clear that
$P_{\bar{a}}$ is a prime ideal which contains $I^p$ and
$\beta(a)=\bar{a}$, where

\[
 a(i)=
\begin{cases}
t_k-1,\quad\text{if}\quad i=j_k \quad \text{for some $k$,}\\
\infty, \quad\text{otherwise.}
\end{cases}
\]

 Now $\bar{a}\in {F}^p$ if and only if
$P_{\bar{a}}\in \Min(I^p)$. Assume $P_{\bar{a}}\not\in \Min(I^P)$.
 Then there is a  prime
ideal $Q\in \Min(I^p)$ such that $Q\subset P_{\bar{a}}$. Suppose
 $Q=(x_{i_1e_1},\ldots,x_{i_he_h})$.  Then
 $\{i_1,\ldots,i_h\}\subset \{j_1,\ldots,j_s\}$ and
 $\{e_1,\ldots,e_h\}\subset\{t_1,\ldots,t_s\}$. On the other hand, since $Q$ is a minimal prime ideal
  of $I^p= \Sect_{i=1}^t J_i^p$, there exists an integer $e\in[t]$ such that $Q$ is one of the minimal prime ideals
  of
 $$
 J_e^p=(x_{i_1}^{b_1},\ldots,x_{i_h}^{b_h})^p.
 $$
It follows that  $1\leq e_i\leq b_i$ for $i=1,\ldots,h$.
 Therefore there exists $b\in {F_e}$  with
 \[
 b(i)=
 \begin{cases}
 e_k-1,\quad\text{if}\quad
 i\in\{i_1,\ldots,i_h\},
 \\\infty,\quad\text{otherwise}.
 \end{cases}
 \]
 This implies that  $a< b\leq m_e$, and $\infpt(a)<\infpt(b)=\infpt(m_e)$,
 a contradiction.

\medskip
 Next we show that $\beta$ is injective: let $a,b\in {F}$ and $a\neq b$.
 Then there exists an integer $i$ such that $a(i)\neq
 b(i)$. We have to show $\bar{a}\neq\bar{b}$. We consider
 different cases:

 (i) If $a(i)=0$, and $b(i)\neq 0$,
 then $\bar{b}(i1)=\infty$ and $\bar{a}(i1)=0$.

(ii) If $a(i)=\infty$, and $b(i)=t-1$ where
 $t\neq\infty$, then  $\bar{a}(it)=\infty$  and  $\bar{b}(it)=0$.

(iii) Suppose   $0<t-1=a(i)\neq \infty$. If
 $b(i)=0$, then we have   case (i). If $b(i)=\infty$ then we have
 case (ii). Finally if  $0<s-1=b(i)\neq \infty$, then  $t\neq s$ since $a(i)\neq b(i)$ and hence
 $\bar{a}(it)=0$ and $\bar{b}(it)=\infty$.

 In all cases it follows that $\bar{a}\neq \bar{b}$.

\medskip
 Finally we show that $\beta$ is surjective: let $\bar{a}\in {F}^p\subset\Union_{i=1}^t F_i^p$
  be any facet of $\Gamma^p$. Then there
exists an integer $i\in[t]$ such that $\bar{a}\in F_i^p$.
Therefore $P_{\bar{a}}$ is a minimal prime ideal of
$$J_i^p=(x_{i_1}^{a_1},\ldots,x_{i_k}^{a_k})^p,$$  and hence there
exists $t_i\leq a_i$ such that
$P_{\bar{a}}=(x_{i_1t_1},\ldots,x_{i_kt_k})$. Therefore
 \[
 \bar{a}(ij)=
 \begin{cases}
 0,\quad\text{if}\quad
 i=i_r\quad\text{and}\quad
 j=t_r\quad\text{for some}\quad r\in[k],
  \\ \infty,\quad\text{otherwise}.
 \end{cases}
 \]

 By our definition we have  $\bar{a}=\beta(a)$, where $a\in A$
 with
 \[
 a(i)=
 \begin{cases}
 t_r-1,\quad\text{if}\quad i=i_r\in\{i_1,\ldots,i_k\},
 \\ \infty,
 \quad\text{otherwise}.
 \end{cases}
 \]

 It will be enough to show that  $a\in {F}$. Since $\bar{a}\in F_i^p$ and the restriction of
 $\beta $ to $F_i$ is a bijection from $F_i$ to $F_i^p$, it follows that $a\in F_i$.
 If $a\not\in {F}$, then there exists some $j\neq i$, such that   $a\leq m_j$,
  and  $\infpt(a)<\infpt(m_j)$.  Therefore there exists an element
 $b\in {F_j}$,  such that $b(i)=a(i)$ for all $i$ with  $b(i)\neq\infty$.
 This implies that $a< b$,  and $\infpt(a)<\infpt(b)=\infpt(m_j)$.
 It follows from the definition of the map $\beta$ that $\bar{a}< \bar{b}$,
 and that $P_{\bar{b}}$ is a prime ideal with $I^p\subset P_{\bar{b}}\subsetneqq
 P_{\bar{a}}$, a contradiction.
 \end{proof}

\medskip
Now let $I\subset S$ be a monomial ideal and $I^p\subset T$ be the
polarization of $I$. Furthermore let
$$\pi :T \longrightarrow S,\quad x_{ij} \mapsto x_i.$$ be the
epimorphism which attached to the polarization. Note that
$$\ker(\pi)=(x_{11}-x_{12},\ldots,x_{11}-x_{1r_1},\ldots,x_{n1}-x_{n2},\ldots,x_{n1}-x_{nr_n})$$
where $r_i$ is the number of variables of the form $x_{ij}$ which
are needed for polarization. Set
$$y:=x_{11}-x_{12},\ldots,x_{11}-x_{1r_1},\ldots,x_{n1}-x_{n2},\ldots,x_{n1}-x_{nr_n}, $$
then $y$ is a sequence of linear forms in $T$.

\begin{Proposition}\label{up&dawon} Let $I\subset S$ be a monomial ideal
and $I^p$ be the polarization of $I$.
Assume that
\[
\mathcal{G}:\; I^p=J_0\subset J_1\subset\ldots\subset J_r=T
\]
 is a clean filtration of $I^p$, and that
 \[
 \mathcal{F}:\; I=I_0\subset I_1\subset\ldots\subset I_r=S
 \]
is the specialization of $\mathcal{G}$, that is, $\pi(J_i)=I_i$
for all $i$. Then $\mathcal{F}$ is a pretty clean filtration of
$I$ with $I_k/I_{k-1}\iso S/\pi(Q_k)$, where $Q_k\iso
J_k/J_{k-1}$.
\end{Proposition}

\begin{proof} For each $k\in[r]$ the $S$-module  $I_k/I_{k-1}$
is a cyclic module since $J_k/J_{k-1}$ is cyclic for all $k$. Let
$I_k/I_{k-1}\iso S/L_k$, where $L_k$ is a monomial ideal in $S$.
It is clear that $\pi(Q_k)\subset L_k$. Indeed, $Q_k=J_{k-1}:u_k$,
where $J_k=(J_{k-1},u_k)$ and where $J_k/J_{k-1}\iso T/Q_k$. If
$v\in Q_k$, then $vu_k\in J_{k-1}$. It follows that
$\pi(vu_k)=\pi(v)\pi(u_k)\in\pi(J_{k-1})=I_{k-1}$, and hence
$\pi(v)\in I_{k-1}:\pi(u_k)=L_k$.

We want to show that $\pi(Q_k)=L_k$. $S$ and $T$ are standard
graded with $\deg(x_i)=\deg(x_{ij})=1$ for all $i$ and $j$, and
$\mathcal{G}$ is a graded prime filtration of $I^p$. Therefore
$\mathcal{F}$ is a graded filtration of  $I$, and we have the
following isomorphisms of graded modules $J_i/J_{i-1}\iso
T/Q_i(-a_i)$ and $I_i/I_{i-1}\iso S/L_i(-a_i)$, where
$a_i=\deg(u_i)=\deg(\pi(u_i))$.

The filtrations  $\mathcal{G}$ and $\mathcal{F}$ yield the
following  Hilbert series of $T/I^P$ and $S/I$:
\[
\Hilb(T/I^p)=\sum_{i=1}^r\Hilb(T/Q_i)t^{a_i}\quad
\text{and}\quad \Hilb(S/I)=\sum_{i=1}^r\Hilb(S/L_i)t^{a_i}.
\]
On the other hand since $y$ is a regular sequence of linear forms
on $T/I^p$ and on $T/Q_i$ for each $i\in[r]$,  we have
\begin{eqnarray*}
\Hilb(S/I)&=&(1-t)^l\Hilb(T/I^p)=(1-t)^l\sum_{i=1}^r\Hilb(T/Q_i)t^{a_i}\\
&=&\sum_{i=1}^r(1-t)^l\Hilb(T/Q_i)t^{a_i}=\sum_{i=1}^r\Hilb(S/\pi(Q_i))t^{a_i},
\end{eqnarray*}
where $l=|y|$.

On the other hand, since $\pi(Q_i)\subset L_i$, we have the
coefficientwise inequality $\Hilb(S/L_i)\leq\Hilb(S/\pi(Q_i))$,
and equality holds if and only if $L_i =\pi(Q_i)$. Therefore we
have
$$\Hilb(S/I)=\sum_{i=1}^r\Hilb(S/\pi(Q_i))t^{a_i}\geq\sum\Hilb(S/L_i)t^{a_i}=\Hilb(S/I).$$
It follows that $L_i=\pi(Q_i)$ is a prime ideal for $i=1,\ldots,
r$.

We know that $\Gamma^p$ the multicomplex associated to $I^p$ is
shellable, since $I^p$ is clean. Therefore we may assume that
$\mathcal{G}$ is obtained from a shelling of $\Gamma^p$. By
\cite[Corollary 10.7]{HP} and its proof (or directly from the
definition of shellings of a simplicial complex) it follows that
$\mu(Q_i)\geq\mu(Q_{i-1})$ for all $i\in[r]$, where $\mu(Q_i)$ is
the number of generators of $Q_i$. Since by Corollary \ref{minI^p}
each $Q_i$ is of the form $(x_{i_1t_1},\ldots, x_{i_kt_k})$ with
$i_r\neq i_s$ for $r\neq s$,  it follows that
$\mu(Q_i)=\mu(\pi(Q_i))=\mu(L_i)$. Therefore $\mu(L_i)\geq
\mu(L_{i-1})$ for all $i$. This implies that $\mathcal F$ is a
pretty clean filtration of $S/I$.
\end{proof}

As the main result of this section we have
 \begin{Theorem}\label{Ghomri} Let $I\subset S=K[x_1,\ldots,x_n]$ be a monomial
 ideal and $I^p$ its polarization. Then the following are equivalent
\begin{enumerate}
\item[(a)] $I$ is pretty clean.
\item[(b)] $I^p$ is clean.
\end{enumerate}
\end{Theorem}

\begin{proof} (a)\implies (b): Assume $I$ is pretty clean. Then the multicomplex $\Gamma$ associated with $I$ is shellable.
Let $a_1,\ldots,a_r$ be a shelling of $\Gamma$, and
\[
\mathcal{F}:\; I=I_0\subset I_1\subset\ldots\subset I_r=S
\]
the  pretty clean filtration of $I$ which is obtain from this
shelling, i.e,  $I_i=\Sect_{k=1}^{r-i} I(\Gamma(a_k))$.  Let
$\mathcal{F}^p$ be the polarization of $\mathcal{F}$. By
Proposition \ref{ffp},  $\mathcal{F}^p$ is a prime filtration of
$I^p$ with $\ell(\mathcal{F})=\ell(\mathcal{F}^p)$. Using
Proposition \ref{martaba} we have
\[
|{F}(\Gamma^p)|=|{F}(\Gamma)|.
\]
On the other hand, since $I$ is pretty clean we know that
$\ell(\mathcal{F})=|{F}(\Gamma)|$. Hence we conclude that
\[
|{F}(\Gamma^p)|=\ell(\mathcal{F}^p).
\]
Therefore,  since $\Min(I^p)=\Ass(I^p)\subset \Supp({\mathcal
F}^p)$, it follows that $\Min(I^p)= \Supp({\mathcal F}^p)$, which
implies  that $I^p$ is clean.

(b)\implies (a) follows from  Proposition \ref{up&dawon}.
\end{proof}

As an immediate consequence we obtain the following result of
\cite[Corollary 10.7]{HP}:

\begin{Corollary} Let $I\subset S$ be a monomial ideal, and
\[
\mathcal{F}:\; I=I_0\subset I_1\subset\ldots\subset I_r=S
\]
a prime filtration of $S/I$ with $I_j/I_{j-1}\iso S/P_j$. Then the
following are equivalent:
\begin{enumerate}
\item[(a)] $\mathcal{F}$ is a pretty clean filtration of $S/I$.
\item[(b)] $\mu(P_i)\geq\mu(P_{i+1})$ for all $i=0,\ldots,r-1$.
\end{enumerate}
\end{Corollary}

\section{A new characterization of  pretty clean monomial ideals}
Let $R$ be a Noetherian ring, and $M$ a finitely generated
$R$-module. For $P\in \Spec(R)$ the number
$\mult_M(P)=\length(H^0_P(M_P))$ is called the {\em length
multiplicity} of $P$ with respect to $M$. Obviously, one has
$\mult_M(P)>0$ if and only if $P\in\Ass(M)$. Assume now that
$(R,\mm)$ is a local ring. Recall that the {\em arithmetic degree}
of $M$ is defined to be
$$\adeg(M)=\sum_{P\in\Ass(M)}\mult_M(P)\cdot \e(R/P),$$ where
$\e(R/P)$ is the {\em multiplicity} of the associated graded ring
of $R/P$.

First we notice the following

\begin{Lemma}\label{mult} Suppose $R$ is a Noetherian ring, and $M$ a finitely
generated $R$-module. Let
\[
\mathcal{F} \: 0=M_0\subset M_1\subset\ldots\subset M_r=M
\]
be a prime filtration of $M$ with $M_i/M_{i-1}\iso R/P_i$. Then
$$\mult_M(P)\leq |\{i\in[r-1]:\; M_{i+1}/M_i\iso R/P\}|$$ for all
$P\in\Spec(R)$.
\end{Lemma}
\begin{proof} If $P\not\in\Ass(M)$, the assertion is trivial.
So now let $P\in\Ass(M)$.  Localizing at $P$ we may assume that
$P$ is the maximal ideal of $M$.

Now we will prove the assertion by induction on
$\ell(\mathcal{F})$.
 If $\ell(\mathcal{F})$=1, then the assertion is obviously true. Let
$\ell(\mathcal{F})>1$. From the following short exact sequence
\[
0\to M_1\to M\to M/M_1\to 0
\]
we get the following long exact sequence
\[
0 \to H^0_P(M_1) \to H^0_P(M) \to H^0_P(M/M_1) \to \ldots
\]
Therefore $\mult_M(P)=\ell(H^0_P(M))\leq
\ell(H^0_P(M_1))+\ell(H^0_P(M/M_1))$. By induction hypothesis
\[
\mult_{M/M_1}(P)=\ell(H^0_P(M/M_1))\leq
|\{i\in[r-1]\backslash\{1\}:\; M_{i+1}/M_i\iso R/P\}|.
\]
Now consider the following two cases:

(i) If $M_1\iso R/P$, then $\ell(H^0_P(M_1))=1$. Therefore
\[
\mult_M(P)\leq 1+\mult_{M/M_1}(P)\leq |\{i \in [r-1] :\;
M_{i+1}/M_i\iso R/P\}|.
\]
(ii) If $M_1\not\iso R/P$, then $\ell(H^0_P(M_1))=0$. Hence
\[
\mult_M(P) \leq |\{i \in [r-1] :\; M_{i+1}/M_i\iso R/P\}|.
\]
\end{proof}

Let $S=K[x_1,\ldots,x_n]$ be the polynomial ring in $n$ variables
over the field $K$. Let $I\subset S$ be a monomial ideal and
$\Gamma$ be the multicomplex associated to $I$. We denote the
arithmetic degree of $S/I$ by $\adeg(I)$. Since $\e(S/P)=1$ for
all $P\in\Ass(I)$, it follows that
$\adeg(I)=\sum_{P\in\Ass(I)}\mult_I(P)$, where
$\mult_I(P)=\mult_{S/I}(P)$. By \cite[Lemma 3.3]{STV}
$\adeg(I)=|\Std(I)|$, where $\Std(I)$ is the set of standard pairs
with respect to $I$. Also by \cite[Lemma 9.14]{HP}
$|\Std(I)=|F(\Gamma)|$. Since $|F(\Gamma)|=|F(\Gamma^p)|$, see
\ref{martaba}, it follows that $\adeg(I)=\adeg(I^p)$, where $I^p$
is the polarization of $I$ and $\Gamma^p$ the multicomplex
associated to $I^p$.

In this part we want to show that $\adeg(I)$ is a lower bound for
the length of any prime filtration of $S/I$ and the equality holds
if and only if $S/I$ is a pretty clean module.
\begin{Theorem}\label{adeg} Let $I\subset S$ be a monomial ideal and
$\mathcal{F}$ a prime filtration of $I$. One has
\begin{enumerate}
\item $\adeg(I)\leq \ell(\mathcal{F})$;
\item $\ell(\mathcal{F})=\adeg(I)\Leftrightarrow \mathcal{F}$ is a
pretty clean filtration of $I$.
\end{enumerate}
\end{Theorem}
\begin{proof} Part 1 is clear by Lemma \ref{mult}.

One direction of (2) is \cite[Corollary 6.4]{HP}. For the other
direction assume
$\ell(\mathcal{F})=\adeg(I)=|F(\Gamma)|=|F(\Gamma^p)|$. By
Proposition \ref{ffp} $\mathcal{F}^p$ is a prime filtration of
$I^p$ with $\ell(\mathcal{F}^p)=|F(\Gamma^p)|$= the number of
minimal prime ideals of $\Gamma^p$. Therefore $\mathcal{F}^p$ is a
clean filtration of $I^p$, so by Theorem \ref{Ghomri}
$\mathcal{F}$ is a pretty clean filtration of $I$.
\end{proof}

Combining Theorem \ref{adeg} with Theorem \ref{Ghomri} we get
\begin{Corollary} Let $I\subset S$ be a monomial ideal.
Assume $\Gamma$ is the multicomplex associated to $I$ and $I^p$
the polarization of $I$. The following are equivalent:
\begin{enumerate}
\item[(a)] $\Gamma$ is shellable;
\item[(b)] $I$ is pretty clean;
\item[(c)] There exists a prime filtration $\mathcal{F}$ of $I$ with
$\ell(\mathcal{F})=\adeg(I)$;
 \item[(d)] $I^p$ is clean;
 \item[(e)] If $\triangle$ be the simplicial complex associated to
 $I^p$, then $\triangle$ is shellable.
 \end{enumerate}
\end{Corollary}

If $R$ is a Noetherian ring and $M$ a finitely generated
$R$-module with  pretty clean filtration $\mathcal{F}$, then
$\Ass(M)=\Supp(\mathcal{F})$, see \cite[Corollary 3.6]{HP}. The
converse is not true in general as shown  by an example in
\cite{HP}.  The example given there is a cyclic module defined by
a non-monomial ideal.  The following example shows that even in
the monomial case the converse does  not hold in general.

\begin{Example}\label{gut} {\em Let $S=K[a,b,c,d]$ be the polynomial ring
over the field $K$, $I\subset S$ the ideal
\[
 I=(a,b)\cdot (c,d)\cdot (a,c,d)=(abc,abd,acd,ad^2,a^2d,ac^2,a^2c,bcd,bc^2,bd^2)
 \]
 and $M=S/I$. We claim that the module $M=S/I$ is not pretty
 clean,  but that $M$ has a prime filtration
 $\mathcal{F}$ with $\Supp(\mathcal{F})=\Ass(M)$.

 Note that
 $(a,b)\cap(c,d)\cap(a,c,d^2)\cap(a,c^2,d)\cap(a^2,b,c,d^2)\cap(a^2,b,c^2,d)$
 modulo $I$ is an irredundant primary decomposition of (0) in
 $M$.

  We see that $\Ass(M)=\{(a,b),(c,d),(a,c,d),(a,b,c,d)\}$.
 It is clear that
\begin{eqnarray*}
\mathcal{F}:\;  I&=&I_0\subset I_1=(I,ac)\subset
I_2=(I_1,ad)\subset I_3=(I_2,bd)\\&\subset& I_4=(I_3,bc) \subset
I_5=(I_4,a)\subset I_6=(a,b)\subset S
\end{eqnarray*}
 is a prime filtration
of $M$ with $\Supp(\mathcal{F})=\Ass(M)$. Indeed $I_1/I\iso
I_2/I_1\iso S/(a,b,c,d)$, $I_3/I_2\iso I_4/I_3\iso I_6/I_5\iso
S/(a,c,d)$ and $I_5/I_4\iso S/(c,d)$.

From the above irredundant primary decomposition of $I$ it follows
that $\adeg(I)$=6. But the length of any prime filtration of $I$
is at least  7. Therefore $I$ can not be pretty clean. In other
words, from \cite[Corollary 1.2]{HP} it follows that
$D_1(M)=((a,b)\cap(c,d))/I$ and that $D_2(M)=M$, where $D_i(M)$ is
the largest submodule of $M$ with $\dim(M)\leq i$, for
$i=0,\ldots,\dim(M)$. It follows that $D_2(M)/D_1(M)\iso
S/(a,b)\cap(c,d)$  is not clean. Knowing now $D_2(M)/D_1(M)$ is
not clean, we conclude from \cite[Corollary 4.2]{HP} that $M=S/I$
is not pretty clean.}
\end{Example}

\end{document}